\documentclass[11pt]{amsart}
\usepackage{amsmath}
\usepackage{amssymb,latexsym}
\usepackage{amscd}
\usepackage{xspace}
\usepackage{verbatim}

\numberwithin{equation}{section}

\newtheorem{sub}{\name}[section]
\newcommand{\bs}{
\begin{sub}}
\newcommand{\es}{
\end{sub}}
\newcommand{\bsl}[1]{
\begin{sub}\label{#1}}
\newcommand{\bth}[1]{\def\name{Theorem}
\begin{sub}\label{t:#1}}
\newcommand{\blemma}[1]{\def\name{Lemma}
\begin{sub}\label{l:#1}}
\newcommand{\bcor}[1]{\def\name{Corollary}
\begin{sub}\label{c:#1}}
\newcommand{\bdef}[1]{\def\name{Definition}
\begin{sub}\label{d:#1}}
\newcommand{\bprop}[1]{\def\name{Proposition}
\begin{sub}\label{p:#1}}

\newcommand{\rth}[1]{Theorem~\ref{t:#1}}

\newcommand{\rdef}[1]{Definition~\ref{d:#1}}

\newcommand{\BA}{\begin{array}}
\newcommand{\EA}{\end{array}}
\newcommand{\BSA}{
\begin{subarray}}
\newcommand{\ESA}{
\end{subarray}}
\newcommand{\BAL}{
\begin{aligned}}
\newcommand{\EAL}{
\end{aligned}}
\newcommand{\BALG}{
\begin{alignat}}
\newcommand{\EALG}{
\end{alignat}}
\newcommand{\BALGN}{
\begin{alignat*}}
\newcommand{\EALGN}{
\end{alignat*}}
\newcommand{\note}[1]{\noindent\textit{#1.}\hspace{2mm}}

\newcommand{\Remark}{\note{Remark}}

\newcommand{\bproof}{\begin{proof}}
\newcommand{\eproof}{\end{proof}}
\newcommand{\forevery}{\quad \forall}

\newcommand{\tin}{\to\infty}

\newcommand{\abs}[1]{\left |#1\right |}

\def\angb<#1>{\langle #1 \rangle}

\newcommand{\rec}[1]{\frac{1}{#1}}
\newcommand{\opname}[1]{\text{\rm #1}\,}

\newcommand{\dist}{\opname{dist}}
\newcommand{\sign}{\opname{sign}}
\newcommand{\myfrac}[2]{{\displaystyle \frac{#1}{#2} }}
\newcommand{\myint}[2]{{\displaystyle \int_{#1}^{#2}}}

\newcommand{\q}{\quad}
\newcommand{\qq}{\qquad}

\newcommand{\1}{\\[1mm]}
\newcommand{\2}{\\[2mm]}
\newcommand{\ity}{\infty}
\newcommand{\prt}{\partial}

\newcommand{\tl}{\tilde}
\newcommand{\sbs}{\subset}

\newcommand{\unl}{\underline}

\newcommand{\Consy}{Consequently\xspace}

\newcommand{\ifif}{if and only if\xspace}
\newcommand{\suff}{sufficiently\xspace}
  \def\ro  {\rho_{_{
\partial\Omega}}}

\def\ga{\alpha}            \def\gg{\gamma}
       \def\gd{\delta}      
                         
\def\gf{\phi}

    \def\gr{\rho}        
\def\gs{\sigma}       
      
                \def\gz{\zeta}
\def\Gg{\Gamma}     \def\Gd{\Delta}

\def\Gw{\Omega}              

   \def\CU{{\mathcal U}}

   \def\BBN {\mathbb N}    
   \def\BBR {\mathbb R}

\newcommand{\bdw}{\prt\Gw\xspace}

\begin{document}
\title[Maximal solutions]{Maximal solutions of semilinear elliptic equations with locally integrable forcing term.}
\author{Moshe Marcus }
\address{Department of Mathematics, Technion\\
 Haifa 32000, ISRAEL}
 \email{marcusm@math.technion.ac.il}
\author{ Laurent Veron}
\address{Laboratoire de Math\'ematiques, Facult\'e des Sciences\\
Parc de Grandmont, 37200 Tours, FRANCE}
\email{laurent.veron@univ-tours.fr}
\thanks{This research was partially supported by an EC Grant through
the RTN Program "Front-Singularities", HPRN-CT-2002-00274.}
\date{\today}
\begin{abstract}
 We study the existence of a maximal solution of $-\Gd u+g(u)=f(x)$
in a domain $\Gw\subset \BBR^N$ with compact boundary, assuming that
$f\in (L^1_{loc}(\Gw))_+$ and that $g$ is nondecreasing, $g(0)\geq
0$ and $g$ satisfies the Keller-Osserman condition. We show that
if the boundary satisfies the classical $C_{1,2}$ Wiener criterion
then the maximal solution is a large solution, i.e., it blows up
everywhere on the boundary. In addition we discuss the question of
uniqueness of large solutions.
\end{abstract}
\maketitle
\section {Introduction}
\setcounter{equation}
{0} Let  $\Gw$ denote a subdomain of  $\mathbb R^N$, $N\geq 2$,
$\ro(x)=\dist(x,\prt\Gw), \forevery x\in \mathbb R^N$, and $g\in C(\mathbb R)$ is
nondecreasing.
In a preceding article \cite {MV3}, we studied existence and
uniqueness of solutions of the problem
\begin {equation}
\label {main-0} -\Gd u+g(u)=0\quad \mbox {in }\;\Gw,
\end {equation}
subject to the boundary blow-up condition
\begin {equation}
\label {blow-up}
\lim_{
\begin {subarray}{c}\ro (x)\to 0\\
x\in K
\end{subarray}}u(x)=\ity\forevery K\subset\Gw,\quad K\mbox { bounded}.
\end {equation}
Such a function $u$ is called a {\em large solution}.
In this article we extend the study to the equation with forcing term,
\begin {equation}
\label {main-f} -\Gd u+g(u)=f(x)\quad \mbox {in }\;\Gw,
\end {equation}
where $f\in L^{1}_{loc}(\Gw)$ is nonnegative. We assume throughout the paper that $g$ satisfies the following
conditions:
\begin{equation}\label{basic}
\begin{aligned}
g\in C(\BBR),\quad g \text{ non decreasing},\quad g(0)\geq 0.
\end{aligned}
\end{equation}
By a solution of (\ref{main-f}) we mean a locally integrable function $u$ such that $g(u)\in
L_{loc}^1(\Gw)$ and (\ref{main-f}) holds in the distribution sense.
Accordingly, if  $u$ is a  solution of (\ref{main-f}) then $\Gd u \in
L^1_{loc}(\Gw)$ and consequently $u\in W^{1,p}_{loc}(\Gw)$ for
some $p>1$ (see~\cite{BS}). Therefore, if  $\Gw'$ is a smooth
bounded domain such that $\bar\Gw'\sbs \Gw$, then
 $u$ possesses an $L^1$ trace on $\bdw'$ and, if $\gf$ is a non-negative function in
 $C_0^2(\bar\Gw')$, i.e., $\gf\in C^2(\bar\Gw')$
 and $\gf=0$ on $\bdw'$, then
\begin{equation}\label{locsol}
 \int_{\Gw'}(-u\Gd\gf+g(u)\gf)\,dx=\int_{\Gw'}f\gf \,dx -\int_{\bdw'}u\prt{\gf}/\prt\mathbf{n'}dS,
\end{equation}
where $\mathbf{n'}$ denotes the external unit normal on $\bdw'$.
The boundary blow-up condition
should be understood as an essential limit: $u$ is bounded below a.e. by a function
$u_0$ which satisfies (\ref{blow-up}).
\par  In a well known paper \cite {Br}
Brezis proved that, for any $q>1$ and $f\in L^1_{loc}(\BBR^N)$,
there exists a unique solution of the equation
\begin {equation}
\label {u^q} -\Gd u+\abs u^{q-1}u=f\quad \mbox {in }\;\BBR^N.
\end {equation}
The proof  was based upon a duality argument which implied local
$L^q_{loc}(\BBR^N)$-bounds of approximate solutions.

In the present paper we  investigate this problem, for $f\geq 0$, for a
large family of nonlinearities and arbitrary domains with compact
boundary satisfying a mild regularity assumption. When $\Gw\subsetneqq
\BBR^N$, we shall concentrate on the existence and uniqueness of {\em large solutions},
i.e. solutions which blow up on the boundary. Other boundary value problems
may have no solution when $f\in (L^1_{loc})_+$. For instance, if $\Gw$ is a smooth, bounded
domain and the boundary data is in $L^1(\bdw)$ then the
boundary value problem for (\ref{main-f}) possesses a solution (in
the $L^1$ sense) \ifif $f\in L^1(\Gw;\gr)$, where
$\gr(x)=\dist(x,\bdw)$. In fact, in this case, if $f\in C(\Gw)$
and $f\geq c_0\gr^{-2}$ for some positive constant $c_0$,
then every solution $u$ of (\ref{main-f}), such that $u\geq 0$ in a
neighborhood of the boundary, is necessarily a large solution. However one can establish a partial
result, namely,
the existence of a {\em minimal solution} of the equation which is also a supersolution of the boundary value
problem, (see \rth{0} below).

\par The problem of
existence of large solutions is closely related to the question of
existence of {\em maximal solutions}. A maximal solution (if it
exists) need not be a large solution. It is well known that, for
equation (\ref{u^q}) with $f=0$,
 a maximal solution exists in any domain.
This is a consequence of the estimates of Keller~\cite{Ke} and Osserman~\cite{Os} as it was shown
in \cite{LN}.
In a recent paper, Labutin~\cite{La}
presented a necessary and sufficient condition on $\Gw$, for the maximal solution
of (\ref{u^q}) with $f=0$ to be a large solution.

 A function $g$ satisfies the Keller-Osserman condition (see \cite{Ke} and \cite{Os})
 if for every $a>0$
\begin {equation}\label {KO}
\int_{a}^\infty\left(\int_{0}^t g(s)\,ds\right)^{-1/2}dt<\infty.
\end {equation}

\par Our first result concerns the existence of maximal solutions.

\bth{I}  Let $\Gw$ be a domain in
 $\mathbb R^N$ and  let $g$ be a function satisfying  (\ref{basic}) and the Keller-Osserman condition.
In addition assume that (\ref {main-0}) possesses a
subsolution. Then  (\ref {main-f}) possesses a maximal solution, for
every  non-negative $f\in L^1_{loc}(\Gw)$.
 \es

\Remark  If $\Gw$ is  bounded  or if $g (r_{0})=0$ for some
$r_{0}\in \BBR$ then equation (\ref{main-0}) possesses a solution. In fact
it possesses a {\em bounded} solution.
\par If $g$ remains positive and the domain is unbounded, some
conditions for the existence of a solution of (\ref{main-0}) can
be found in \cite {MV3}.

\par The existence of a maximal solution implies that the family of all solutions
of (\ref{main-f}) is locally uniformly bounded from above.
By \cite{Ke} and \cite{Os} the Keller Osserman condition is {\em necessary} for this property
to hold. Furthermore this property implies that a family of solutions which is locally uniformly
bounded from below is compact.

\par In the next result we consider boundary value problems with $L^1$ boundary data.
\bth{0} Suppose that $g$ satisfies (\ref{basic}) and the Keller-Osserman condition.\\
(i) Assume that $\Gw$ is a  smooth bounded domain, $f\in (L^1_{loc})_+(\Gw)$ and  $h\in L^1(\bdw)$.
 Then there
exists a {\em minimal supersolution} $\unl{u}_h\in L^1_{loc}(\Gw)$ of the boundary value problem
\begin{equation}\label{Lbvp}
 -\Gd u+g(u)=f \q\text{in}\q\Gw, \qq u=h\q\text{on}\q\bdw.
\end{equation}
The function $\unl{u}_h$ satisfies (\ref{main-f}) and, if $f\in L^1(\Gw;\gr)$, it
is the unique solution of (\ref{Lbvp}).\\
(ii) Assume that $\Gw$ is a  bounded domain satisfying the classical Wiener condition,
$f\in (L^1_{loc})_+(\Gw)$ and  $h\in C(\bdw)$.  Then there
exists a {\em minimal supersolution} $\unl{u}_h\in L^1_{loc}(\Gw)$ of (\ref{Lbvp}).
The function $\unl{u}_h$ satisfies (\ref{main-f}) and, if $f\in L^\infty(\Gw)$, it
is the unique solution of (\ref{Lbvp}).
\es
For the definition of a supersolution of the boundary value problem (\ref{Lbvp}) when $f$
is only locally integrable see Section 3. The definition of a sub/super solution of  {\em equation}
(\ref{main-f}) is standard:
\bdef{sub-super}{\rm A function $u\in L^1_{loc}(\Gw)$ is a  subsolution
(resp. supersolution)
of {\em equation} (\ref{main-f}), with $f\in L^1_{loc}(\Gw)$,  if $g(u)\in L^1_{loc}(\Gw)$ and
$$ -\Gd u +g(u)-f\leq 0 \q(\text{resp. }\geq 0) \q\text{in }\Gw$$
in the distribution sense.}
\es
\par We note that if $u$ is a  supersolution of equation (\ref{main-f}), there exists a positive Radon
measure $\mu$ in $\Gw$ such that
$$  -\Gd u +g(u)-f=\mu, \q\text{in}\q\Gw.$$
Therefore (\ref{locsol}) holds with $f$ replaced by $f+\mu$:
$$ \int_{\Gw'}(-u\Gd\gf+(g(u)-f)\gf)\,dx=\int_{\Gw'}\gf \,d\mu -\int_{\bdw'}u\prt{\gf}/\prt\mathbf{n'}dS.$$
 \par The following result concerns the existence of large solutions.
\bth{II}  Let $\Gw$ be a domain in
 $\mathbb R^N$ with  non-empty, compact boundary.
 Assume that $g$ satisfies  (\ref{basic}) and the Keller-Osserman condition and that
(\ref {main-0}) possesses a subsolution $V$ in $\Gw$. Put
$$\CU_V(\Gw):=\{h\in L^1_{loc}(\Gw):\, h\geq V \;a.e.\}.$$
Under these assumptions:\1
(i) For every $f\in (L_{loc}^1)_+(\Gw)$,  (\ref{main-f}) possesses a minimal solution $V_f$ in $\CU_V(\Gw)$.
 $V_f$ increases as $f$ increases.\1
(ii) Assume, in addition, that $\Gw$ satisfies the (classical) Wiener
 criterion. Then, for every $f\in (L_{loc}^1)_+(\Gw)$, (\ref{main-f}) possesses a large solution.
Moreover there exists a {\em minimal large solution} of
(\ref{main-f}) in $\CU_V(\Gw)$.\1
 (iii) If $\Gw$ is bounded and
satisfies the (classical) Wiener
 criterion then (\ref{main-f}) possesses a {\em minimal large solution}.
\es
\Remark (a) Part (i) implies that if (\ref{main-0}) possesses a large solution then
 (\ref{main-f}) possesses a large solution for every $f\in (L_{loc}^1)_+(\Gw)$. In \cite{MV3} it was shown that,
 if $g$ satisfies (\ref{basic}) and the so called
{\em weak singularity condition} then (\ref{main-0}) possesses a large solution  in any domain
$\Gw$ such that $\bdw=\prt\bar\Gw^c$. The weak singularity condition is satisfied, for example,
in the following cases:\\[1mm]
(1)  If $g(u)=|u|^{q-1}u$ and $1<q<N/(N-2)$ for $N\geq3$. \\[1mm]
(2) If $0<g(u)<ce^{au}$, $a>0$, for $N=2$.\\[1mm]
(b) Labutin \cite{La} studied power nonlinearities, $g(u)=|u|^{q-1}u$, $q>1$, and showed that a necessary
and sufficient condition for the existence of large solutions of (\ref{main-0}) is that $\Gw$ satisfy
a Wiener type condition in which the classical capacity $C_{1,2}$ is replaced by the capacity
$C_{2,q'}$. Labutin's condition is less restrictive than the classical Wiener condition; however the
latter applies to every nonlinearity satisfying the conditions of \rth{II}.
It is interesting to know if the classical Wiener condition is necessary
for the existence of large solutions under these general conditions. More precisely we ask: \\[1mm]
\note{Open problem 1} Let $\Gw$ be a bounded domain which does not satisfy
 the (classical) Wiener criterion
at some point $P\in \bdw$. Does there exist a function $g$ satisfying (\ref{basic}) and the Keller-Osserman condition
such that the maximal solution of (\ref{main-0}) is not a large solution?\\[1mm]
\indent In continuation we consider the question of uniqueness of large solutions,
for nonlinearities $g$ as in \rth{II}.
In order to deal with this question in possibly unbounded domains
we have to restrict ourselves to  solutions which are essentially bounded below by a subsolution of
(\ref{main-0}).
\bth{III}  Let $\Gw$ be a domain in
 $\mathbb R^N$ with  non-empty, compact boundary. Assume that $g$ is convex and satisfies  (\ref{basic})
and the Keller-Osserman condition.\1
(i) Let $V$ be a  subsolution of
(\ref {main-0}). If (\ref {main-0}) possesses a unique  large solution in $\CU_V(\Gw)$
 then, for every
$f\in (L_{loc}^1)_+(\Gw)$, (\ref{main-f}) possesses a unique large solution in $\CU_V(\Gw)$.\1
(ii) Let  $\Gw$ be a bounded domain. If
(\ref {main-0}) possesses a unique  large solution then, for every
$f\in (L_{loc}^1)_+(\Gw)$, (\ref{main-f}) possesses a unique large solution.
\es
\Remark Assertion (ii) implies that if (\ref{main-0}) possesses a unique large solution $W$, then (\ref{main-f}) possesses a
unique large solution
bounded below by $W$. However, if $\Gw$ is unbounded,  (\ref{main-f}) may possess additional
large solutions which are not bounded below by $W$.\\[1mm]
Combining the above result with  \cite[Theorem 0.3]{MV3} we obtain the following.
\bcor{III'}
Let $\Gw$ be a bounded domain in
$\mathbb R^N$ such that $\partial\Gw$ is a locally continuous graph. Suppose that $g$ is convex
and satisfies  (\ref{basic}), the Keller-Osserman condition and the superaddivity condition:
\begin {equation}\label {supadd}
g(a+b)\geq g(a)+g(b)-L\,,\forevery a,\,b\geq 0,
\end{equation}
for some $L>0$.
\par Under these conditions,   (\ref{main-f}) possesses at most one
large solution, for every $f\in (L_{loc}^1)_+(\Gw)$.
\par If, in addition, $\bdw$ is bounded then
(\ref{main-f}) possesses exactly one large solution, for every $f$ as above.
\es
\par Finally we present two results involving solutions in the whole space $\BBR^N$.
\bth{IV} Let $\Gw=\BBR^N$. Assume that $g$  satisfies  (\ref{basic})
and the Keller-Osserman condition and that (\ref {main-0}) possesses a subsolution
$V$. Then:\\[1mm]
(i) For every
$f\in (L_{loc}^1)_+(\BBR^N)$, (\ref{main-f}) possesses a solution $u$ in $\CU_V(\BBR^N)$.\\[1mm]
(ii) Assume, in addition, that $g$ is convex.
 If (\ref {main-0}) possesses a {\em unique  solution} in $\CU_V(\BBR^N)$ then, for every
$f\in (L_{loc}^1)_+(\Gw)$, (\ref{main-f}) possesses a {\em unique solution} in $\CU_V(\BBR^N)$.
\es
For the statement of the next theorem we need the following notation.
If $g$ is a function defined on $\BBR$ such that $g(0)=0$,  we
denote by $\tilde g$ the function given by $\tilde g(t)=-g(-t)$ for every real $t$.
\bth{V} Assume $\Gw=\BBR^N$. Suppose that $g$ and $\tilde g$ satisfy (\ref{basic}) and the Keller-Osserman condition.
Then:\\[1mm]
(i) For every $f\in L_{loc}^1(\BBR^N)$,  (\ref{main-f}) possesses a solution.\\[1mm]
(ii) Assume, in addition, that $g$ is convex in $(0,\infty)$ and $g(0)=0$. Then,
for every
$f\in (L_{loc}^1)_+(\BBR^N)$, (\ref{main-f}) possesses a unique positive solution.
\es
\Remark It can be shown that if, in addition to the assumptions of part (ii),
$g$ satisfies the  condition
\begin{equation}\label{s-additivity}
 \rec{c}g(a+b)\leq g(a)+g(b)\leq cg(a+b) \forevery a,b\in (0,\infty)
\end{equation}
for some constant $c>0$, then (\ref{main-f}) possesses a unique solution in $\BBR^N$,
for every $f\in L_{loc}^1(\Gw)$. This condition
means that $g$ behaves essentially like a power. In the case of powers this result is due to Brezis \cite{Br}.\\[1mm]
\note{Open problem 2} For $\ga>0$, let $g_\ga$ be given by
$$g_\ga(t)=(e^{(t^\ga)}-1)\sign t \forevery t\in \BBR.$$
Does there exist $\ga>0$ such that (\ref{main-f}), with $g=g_\ga$,
 possesses a unique solution in $\BBR^N$, for every $f\in L_{loc}^1(\BBR^N)$ ?

 \section {Existence of a maximal solution}
 \bproof[Proof of \rth{I}]  Let $\{\Gw_{n}\}$ be a sequence of
 bounded subsets of $\Gw$ with smooth boundary such that
\begin{equation}\label{exhaustion}
  \Gw_n\uparrow \Gw,\qq \bar\Gw_n\sbs \Gw_{n+1}.
\end{equation}
 For every $n\in \BBN$ and $m,k>0$ denote by $u_{n,m,k}$  the classical solution of
 \begin {equation}\label {approx1}
 -\Gd u+g(u)=f_k:=\min(f,k) \quad\text{in }\;\Gw_{n},\qquad
u=m \quad\text{on }\;\prt\Gw_{n}.
\end {equation}
Further denote by $v_{n,m}$ and $w_{n,k}$ the solutions of
\begin{equation}\label{app1}
 -\Gd v+g(v)=0\quad\text{in }\;\Gw_{n},\qquad
v=m \quad\text{on }\;\prt\Gw_{n},
\end{equation}
and
\begin{equation}\label{app1'}
-\Gd w=f_k\quad\text{in }\;\Gw_{n},\qquad
w=0 \quad\text{on }\;\prt\Gw_{n},
\end{equation}
respectively.
Then $u_{n,m,k}-v_{n,m}\geq 0$ and hence
$$-\Gd (u_{n,m,k}-v_{n,m})=f_k-g(u_{n,m,k})+g(v_{n,m})\leq f_k.$$
Since $u_{n,m,k}-v_{n,m}$ vanishes on $\prt\Gw_n$, it follows that
\begin{equation}\label{app2}
 u_{n,m,k}- v_{n,m}\leq w_{n,k} \forevery m\in \BBN.
\end{equation}
Both $m\mapsto v_{n,m}$ and $m\mapsto u_{n,m,k}$ are increasing and $v_{n,m}\leq u_{n,m,k}$.
If $g$ satisfies the Keller-Osserman condition then $\lim_{m\to\infty} v_{n,m}=v_n$ is the minimal
 large solution of (\ref{main-0}) in $\Gw_n$. Therefore, by (\ref{app2}),
\begin{equation}\label{app2'}
 v_n\leq u_{n,k}=\lim_{m\to\infty} u_{n,m,k}\leq v_n+w_{n,k}.
\end{equation}
Since $w_{n,k}$ is bounded and $v_n$ is locally bounded it follows that $u_{n,k}$ is locally bounded
in $\Gw_n$. Thus $u_{n,k}$ is a large solution of (\ref{approx1}), for every $k>0$.
Both $k\mapsto u_{n,k}$ and $k\mapsto w_{n,k}$ are increasing.
Hence, letting $k\to\infty$ we obtain,
\begin{equation}\label{app2''}
v_n\leq u_n=\lim_{k\to\infty}u_{n,k}\leq v_n+w_{n},
\end{equation}
where $w_n$ is the solution of
\begin{equation}\label{app1''}
-\Gd w=f\quad\text{in }\;\Gw_{n},\qquad
w=0 \quad\text{on }\;\prt\Gw_{n}.
\end{equation}
 For every $\gz\in C^2_c(\Gw_n)$,
$$\int_{\Gw_n}(-u_{n,k}\Gd\gz+g(u_{n,k})\gz)\,dx=\int_{\Gw_n}f_k\gz \,dx.$$
Since $g(u_{n,k})\uparrow g(u_n)$, $f\in L^1(\Gw_n)$ and, by (\ref{app2''}) and (\ref{app1''}),
$u_n\in L^1_{loc}(\Gw_n)$, it follows that,
$$\int_{\Gw_n}(-u_{n}\Gd\gz+g(u_{n})\gz)\,dx=\int_{\Gw_n}f\gz \,dx,$$
for every $\gz\in C^2_c(\Gw_n)$, $\gz\geq 0$. In addition,  $u_n\geq v_n$ and consequently
the negative part of $u_n$ is  bounded. Therefore, if $\Gw_n^+=\Gw_n\cap \{u_n\geq 0\}$, we obtain
$$0\leq \int_{\Gw_n^+}g(u_{n})\gz\,dx<\infty,$$
for every $\gz$ as above. This implies that $g(u_n)\in L^1_{loc}(\Gw_n)$ and $u_n$ is a large
solution of (\ref{main-f}) in $\Gw_n$.
 Clearly $\{u_n\}$ is monotone
decreasing and $u_n\geq v_0$ in $\Gw_n$ for any subsolution $v_0$ of (\ref{main-0});
by assumption such a subsolution exists. Therefore $u:=\lim u_n$ is
well defined and it is a solution of (\ref{main-f}) in $\Gw$. In fact
$u$  is the maximal solution of (\ref{main-f}) in
$\Gw$. Indeed, if  $U$ is a  solution of (\ref{main-f}) then, in view of (\ref{locsol}),
 $U\leq u_n$ in $\Gw_n$, so that $U\leq u$.
 \eproof
\section{Minimal supersolutions of boundary value problems}
We start with the definition of a supersolution of (\ref{Lbvp})
 when $f$ is only locally integrable.
\bdef{supsol}{\rm
Under the conditions of part (i) (resp. (ii)) of \rth{0},
a function $u\in L^1_{loc}(\Gw)$ is a supersolution of the {\em boundary value problem}
(\ref{Lbvp}) if it is
a supersolution of (\ref{main-f}) and,
for every $f_0\in L^1_+(\Gw)$ (resp. $f_0\in L_+^\infty(\Gw)$) such that $f_0\leq f$, $u$ dominates the solution of the
boundary value problem}
$$  -\Gd u+g(u)=f_0 \q\text{in}\q\Gw, \qq u=h\q\text{on}\q\bdw.$$
\es

\bproof[Proof of \rth{0}] First we verify the following assertion:
\par If $u\in L^1_{loc}(\Gw)$ is a supersolution (in the sense of
\rdef{supsol}) of the boundary value problems
\begin{equation}\label{Lbvp-k}
 -\Gd u+g(u)=f_k=\min(f,k) \q\text{in}\q\Gw, \qq u=h\q\text{on}\q\bdw,
\end{equation}
for every $k>0$, then $u$ is a supersolution of (\ref{Lbvp}).
 \par Under the assumptions of part (ii) the assertion is true by definition. Therefore we assume the
 conditions of part (i). Let $\tl f\in L^1_+(\Gw)$ be a function
dominated by $f$ and put $\tl f_{k}:=\min(\tl f,k)$. If  $\tl u_{k}$ is the solution
of (\ref{Lbvp-k}) with $f_k$ replaced by $\tl f_k$ then  $\tl u_{k}\uparrow \tl u$ where
$\tl u$ is the solution of
$$ -\Gd u+g(u)=\tl f \q\text{in}\q\Gw, \qq u=h\q\text{on}\q\bdw.$$
By assumption, $\tl u_k\leq u$, for every $k>0$.  Hence $\tl u\leq u$ and the assertion is proved.
\par Denote by $u_k$ the unique solution of (\ref{Lbvp-k}).
Since $\Gw$ is bounded there exists a solution of (\ref{main-0}). Therefore, by \rth{I},
there exists a maximal solution $\bar u_f$ of (\ref{main-f}).
Then $u_k\leq \bar u_f$ and $\{u_k\}$ is increasing. Consequently $u=\lim u_k$ is a solution of (\ref{main-f})
and by the first part of the proof it is a supersolution of (\ref{Lbvp}). Obviously it is
the minimal supersolution of (\ref{Lbvp}).
\eproof
 \section {Existence of a large solution}
 We recall that an open subset $\Gw$ of $\BBR^{N}$ satisfies the
 Wiener criterion if, for every $\gs\in \prt\Gw$,
  \begin {equation}
\label {wiener}
  \int_{0}^{1}C_{1,2}(B_{s}(\gs)\cap \Gw^c)\myfrac {ds}{s^{N-1}}=\ity,
  \end {equation}
where $C_{1,2}$ stands for the classical (electrostatic) capacity. If $\Gw$ is a domain with compact,
non-empty boundary and the Wiener criterion
 is fulfilled, then for any $\phi\in C(\prt\Gw)$ and $\psi\in
L^\ity_{loc}(\overline\Gw)$, a weak solution of
 \begin {equation}\label {wiener1}\left\{
 \begin{aligned}
 -\Gd w&=\psi\quad\text{in }\;\Gw\\[2mm]
w&=\phi\quad\text{on }\;\prt\Gw,
 \end{aligned}\right.
\end {equation}
is continuous up to $\prt\Gw$.
\par Suppose that  $V$ is a subsolution of (\ref{main-0}), i.e., $V$ and $g(V)$ are in $L^1_{loc}(\Gw)$ and
$-\Gd V +g(V)$ is a negative distribution. It follows that there exists a positive Radon measure
$\mu$ such that
$$-\Gd V+g(V)=-\mu\q\text{in}\q \Gw.$$
Consequently $V\in W^{1,p}_{loc}(\Gw)$ for some $p>1$ and, if  $\Gw'$ is a smooth
bounded domain such that $\bar\Gw'\sbs \Gw$, then $V$ possesses an $L^1$ trace on $\bdw'$ and
\begin{equation}\label{locsubsol}
 \int_{\Gw'}(-V\Gd\gf+g(V)\gf)\,dx=-\int_{\Gw'}\gf \,d\mu -\int_{\bdw'}V\prt{\gf}/\prt\mathbf{n'}dS,
\end{equation}
for every $\gf\in C_0^2(\bar\Gw')$, where $\mathbf{n'}$ denotes the external unit normal on $\bdw'$.
\bproof[Proof of \rth{II}(i)] Let $V$ be a subsolution of (\ref{main-0}) and let  $\{\Gw_n\}$ be
as in the proof of \rth{I}. Let $V_{f,n}$ be the (unique) solution of the problem
\begin {equation}\label {sub1}
 -\Gd u+g(u)=f \quad\text{in }\;\Gw_{n},\qquad
u=V \quad\text{on }\;\prt\Gw_{n}.
\end {equation}
Since $V$ is a subsolution
\begin{equation}\label{sub2}
V_{f,n+1}\geq V_{f,n}\q \text{in}\q\Gw_n.
\end{equation}
By \rth{I}, there exists a
 maximal solution $\bar u_{f,n}$ (resp. $\bar u_f$)
of (\ref{main-f}) in $\Gw_n$ (resp. $\Gw$).
Clearly
\begin{equation}\label{max-decreases}
\bar u_f\big|_{\Gw_n}\leq \bar u_{f,n+1}\big|_{\Gw_n}\leq \bar u_{f,n}.
\end{equation}
Therefore $\{\bar u_{f,n}\}$ converges and the limit $U$ is a solution in $\Gw$ such that
$U\geq \bar u_f$.
As  $\bar u_f$ is the maximal solution it follows that $U=\bar u_f$; thus
\begin{equation}\label{u_f=limit}
\bar u_f= \lim_{n\to\infty}\bar u_{f,n}.
\end{equation}
Since $V_{f,n}\leq \bar u_{f,n}$, (\ref{sub2}) and (\ref{u_f=limit}) imply that the sequence
 $\{V_{f,n}\}$ converges to a solution $V_f$ of (\ref{main-f}). Clearly $V_f$ is the minimal solution
 in $\CU_V$. By the maximum principle, $V_{f,n}$ increases with $f$. Therefore $V_f$ increases with $f$.\\[2mm]
\note{Proof of \rth{II}(ii)}
Let $\{\Gw_n\}$ be a sequence of domains contained in $\Gw$ satisfying (\ref{exhaustion}), such that,
for each $n\in \BBN$,  $\Gg_n=\bdw_n$ is a smooth compact surface. Note that if $\Gw$ is
unbounded then $\Gw_n$  is also unbounded. In this case, let $\{D_{n,j}:\,n,j\in \BBN\}$ be a family of smooth
 bounded domains such that
$$\bar D_{n,j}\sbs D_{n+1,j+1},\qquad
 \prt D_{n,j}=\Gg_n\cup\Gg'_{j}, \qquad \Gg_{n}\cap\Gg'_{j}=\emptyset,$$
where $\Gg_n$ and $\Gg'_j$ are smooth, compact surfaces and
$$\cup_{j\geq 1} D_{n,j}=\Gw_n.$$
Denote
$$\Gw'_j:=\cup_{n\geq 1} D_{n,j}.$$
If $\Gw$ is bounded we put $D_{n,j}=\Gw_n$, $\Gg'_j=\emptyset$ for every $j\in \BBN$ so that, in this case,
$\Gw'_j=\Gw$.
\par Let $V_0$ be the minimal solution of (\ref{main-0}) bounded below by $V$.
Let $u^0_{m,n,j}$ be the solution of the problem
\begin {equation}
\label {umnj}\left\{ \begin{aligned}
 -\Gd u+g(u)&=0\quad&&\text{in }\;D_{n,j}\\
u&=\max(m,V_0)\quad&&\text{on }\;\Gg_n\\
u&=V_0\quad&&\text{on }\;\Gg'_j.
 \end{aligned}\right.
\end {equation}
By the maximum principle, $u^0_{m,n,j}$ increases with $m$ and $j$ and $u^0_{m,n,j}\geq V_0$. In addition,
by the Keller-Osserman estimate, the set
$$\{u^0_{m,n,j}:\, m\geq 1,\;n>n_0, \;j>j_0\}$$
is bounded in $D_{n_0,j_0}$. Therefore there exists a subsequence $\{n'\}$ such that the limit
\begin{equation}\label{zmj-0}
 z^0_{m,j}=\lim_{n'\tin}u^0_{m,n',j}
\end{equation}
exists in $\Gw'_j$,  $z^0_{m,j}$ is a solution of (\ref{main-0}) in this domain and
\begin{equation}\label{z0mj}
 z^0_{m,j}\geq m\q\text{on }\bdw, \qq z^0_{m,j}=V_0\q\text{on } \Gg'_j,\qq z^0_{m,j}\geq V_0\q\text{in } \Gw'_j.
\end{equation}
In fact, if $w^0_{m,j}$ is the solution of the problem
\begin {equation}
\label {wmj}\left\{ \begin{aligned}
 -\Gd w+g(w)&=0\quad&&\text{in }\;\Gw'_j\\
w&=m \quad&&\text{on }\;\bdw, \\
w&= V_0 \quad&&\text{on }\;\Gg'_j,
 \end{aligned}\right.
\end {equation}
then  $w^0_{m,j}\in C(\bar\Gw'_j)$ (Here we use the fact that $\Gw$ satisfies the Wiener criterion.)
In addition, for any $\gd>0$, if $n$ is \suff large then $\Gg_n$ is contained in a $\gd$-neighborhood of $\bdw$. Therefore
$\sup w^0_{m,j}\big|_{\Gg_n}\to m$ as $n\tin$ and $u^0_{m,n,j}\geq w^0_{m,j}$ for all \suff large $n$.
\Consy  $z^0_{m,j}\geq w^0_{m,j}$.
Further, if $U$ is a large solution of (\ref{main-0}) and $U\geq V_0$ then $U$ dominates $u^0_{m,n,j}$ for all \suff large $n$.
Hence $U\geq z^0_{m,j}$.
Therefore
$$\unl u^0_{_V}:=\lim_{j\tin}\lim_{m\tin} z^0_{m,j}$$
is the minimal large solution of (\ref{main-0}) which dominates $V_0$ (and hence $V$).
Consequently,
if $ \unl u^f_{_V}$ denotes the minimal  solution of
(\ref{main-f}) which dominates  $\unl u^0_{_V}$ then $ \unl u^f_{_V}$ is a large solution
of (\ref{main-f}) which dominates $V$. Further, if $U^f\in \CU_V$ is a large solution of
(\ref{main-f}) then, for fixed
$m,j\in \BBN$, $U^f\geq u^0_{m,n,j}$ for all sufficiently large $n$. Hence $U^f\geq z^0_{m,j}$,
which in turn implies $U^f\geq \unl u^0_{_V}$ and $U^f\geq \unl u^f_V$.
Thus $\unl u^f_V$ is the minimal large solution of (\ref{main-f}) in $\CU_V$.

\par For later reference we observe that, for an appropriate choice of $\{D_{n,j}\}$,
\begin{equation}\label{minlarge-f}
 \unl u^f_{_{V}}=\lim_{j\tin}\lim_{m\tin}\lim_{n\tin} u^f_{m,n,j}.
\end{equation}
Of course the family of domains $\{D_{n,j}\}$ can be chosen so that (\ref{minlarge-f}) holds for a
given finite set of functions $f$.\1
\note{Proof of \rth{II}(iii)}
 Put $f_k:=\min(f,k)$, $k\in \BBN$.
 Let $u_{k,m}$ be the (unique) solution of the problem,
 \begin {equation}\label{data=m}
  \begin{aligned}
 -\Gd w+g(w)&=f_k\quad&&\text{in }\;\Gw\\
w&=m\quad&&\text{on }\;\bdw.
 \end{aligned}
\end {equation}
Obviously, $u_{k,m}\leq \bar u_f$ (=the maximal solution of (\ref{main-f})).
Since $m\mapsto u_{k,m}$ is increasing it follows that  $u_k:=\lim_{m\to\infty} u_{k,m}\leq \bar u_f$
is a large solution of $ -\Gd w+g(w)=f_k$ in $\Gw$.
Further,
$k\mapsto u_k$ is also increasing. Thus $\underline{u}^f:=\lim u_k$ is a large solution
of (\ref{main-f}). Every large solution $U$ of (\ref{main-f})
dominates $u_{k,m}$. Therefore $\underline{u}^f$ is the minimal large solution.
\eproof

\section{Uniqueness}
\bproof[Proof of \rth{III}(i)] Let $\{D_{n,j}\}$  be as in the proof of \rth{II}, chosen so that (\ref{minlarge-f})
holds for both $f$ and the zero function. In fact we shall use all the notation introduced in this
proof. Let $U^f_{m,n,j}$ be the solution of the problem
\begin {equation}
\label {Umnj-f}\left\{ \begin{aligned}
 -\Gd u+g(u)&=f\quad&&\text{in }\;D_{n,j}\\
w&=\max(m,V_0)\quad&&\text{on }\;\prt D_{n,j}.
 \end{aligned}\right.
\end {equation}
Then  $U^f_{n,j}=\lim_{m\tin}U^f_{m,n,j}$ is a large solution of (\ref{main-f}) in $D_{n,j}$ and
\begin{equation}\label{maxsol-f}
\bar u_f:=\lim_{j\tin}\lim_{n\tin}\lim_{m\tin}U^f_{m,n,j}
\end{equation}
is the maximal solution of (\ref{main-f}) in $\Gw$. If  (\ref{main-0}) possesses a large solution
then, of course, $\bar u_f$ (resp. $\bar u_0$) is the maximal large solution of (\ref{main-f})
(resp. (\ref{main-0})).
\par Put
$$Z^f=Z^f_{m,n,j}:= U^f_{m,n,j}- u^f_{m,n,j}\leq 0.$$
Then
$$\Gd (Z^{f}-Z^{0})=\,g(U^f_{m,n,j})-g(U^0_{m,n,j})-g(u^f_{m,n,j})+g(u^0_{m,n,j}),$$
in $D_{n,j}$.
We rewrite the right hand side in the form
$$
\bar d_f (U^f_{m,n,j}-U^0_{m,n,j})-\underline d_f (u^f_{m,n,j}-u^0_{m,n,j}),
$$
where
$$
\bar d_f=\myfrac
{g(U^f_{m,n,j})-g(U^0_{m,n,j})}{U^f_{m,n,j}-U^0_{m,n,j}}
,\quad  \underline d_f \myfrac{g(u^f_{m,n,j})-g(u^0_{m,n,j})}{u^f_{m,n,j}-u^0_{m,n,j}}.
$$
 Since $g$ is convex and nondecreasing,
$$
\bar d_f\geq  \underline d_f\geq 0, \qquad\Gd (Z^{f}-Z^{0})\geq \underline d_{f}(Z^{f}-Z^{0})
$$
in $D_{n,j}$.
As $Z^{f}-Z^{0}$=0 on $\prt D_{n,j}$, it follows that
$$Z^{f}- Z^{0}\leq 0\quad\text{in } D_{n,j}.$$
Thus
$$ U^f_{m,n,j}- U^0_{m,n,j}\leq u^f_{m,n,j}- u^0_{m,n,j}$$
and consequently,
$$ \lim_{j\tin}\lim_{n\tin}\lim_{m\tin}(U^f_{m,n,j}- U^0_{m,n,j})\leq
\lim_{j\tin}\lim_{m\tin}\lim_{n\tin}(u^f_{m,n,j}- u^0_{m,n,j}).$$
Hence, by  (\ref{maxsol-f}) and (\ref{minlarge-f}):
$$\bar u_f- \bar u_0\leq \underline u^f_{_{V}}- \underline u^0_{_{V}}.$$
Thus
\begin{equation}\label{Zf<Z0}
 0\leq \bar u_f- \underline u^f_{_{V}}\leq \bar u_0- \underline u^0_{_{V}}.
\end{equation}
Assuming that (\ref{main-0}) possesses a unique large solution dominating $V$, we find that
$\bar u_0- \underline u^0_{_{V}}$ and hence
$\bar u_f= \underline u^f_{_{V}}.$
Therefore (\ref{main-f}) possesses a unique  large solution in the class of functions dominating $V$.\2
\note{Proof of \rth{III}(ii)}
If (\ref{main-0}) possesses a large solution $U_0$ then (\ref{main-f}) possesses
a large solution $U\geq U_0$. Since $\Gw$ is bounded, (\ref{main-f}) possesses a minimal large solution $\unl u^f$
(by \rth{II}(iii)) and a maximal solution $\bar u_f$ (by \rth{I}). If $U_0$ is the unique large solution of (\ref{main-0})
then, by the same argument as in part (i), $\unl u^f=\bar u_f$.
\eproof
\section{Solutions in the whole space}
\bproof[Proof of \rth{IV}]
(i) Let $u^f_R$ be the maximal solution of (\ref{main-f}) in $B_R=B_R(0)$;
its existence is guaranteed by \rth{I}.
If $V$ is a subsolution of (\ref{main-0}),
 $u^f_R\geq V$ and $u_R$ decreases
 with $R$. Hence $u^f=\lim_{R\tin} u^f_R$ is a solution of
(\ref{main-f}) in $\BBR^N$ and $u^f\geq V$.\\[1mm]
(ii) Obviously, $u^f$ is the maximal solution of (\ref{main-f}) in $\BBR^N$.
Next we construct the  minimal solution bounded below by $V$. For $R>0$, let $v^f_R$ be the
solution of the problem
\begin {equation}
\label {ukR}\left\{ \begin{aligned}
 -\Gd v+g(v)&=f\quad&&\text{in }\;B_R\\
v&=V\quad&&\text{on }\;\prt B_R.
 \end{aligned}\right.
\end {equation}
Then
\begin{equation}\label{B-ineq1}
   V\leq v^f_R\leq u^f_R.
\end{equation}
 Since $V$ is a subsolution, $v^f_R$ increases with $R$. Therefore
\begin{equation}\label{RN-ineq1}
  V\leq v^f:=\lim_{R\tin}v^f_R\leq u^f.
\end{equation}
Clearly $v^f$ is the  minimal solution of (\ref{main-f}) bounded below by $V$.

\par As in the proof of \rth{III} we obtain,
$$u^f-v^f\leq u^0-v^0.$$
If (\ref{main-0}) possesses a unique solution in $\BBR^N$ then $u^0=v^0$ and
consequently $u^f=v^f$. Thus (\ref{main-f}) possesses a unique solution bounded below by $V$.
\eproof
\bproof[Proof of \rth{V}]\ \\[1mm]
\note{A-priori Estimates} If $u$ is a solution of (\ref{main-f}) in $\BBR^N$ then
$\tilde u(\cdot)=-u(-\cdot)$ satisfies
\begin{equation}\label{modified-g}
-\Gd \tilde u +\tilde g(\tilde u)=\tilde f \quad\text{in }\BBR^N,
\end{equation}
where $\tilde f(x)=-f(-x)$.
\par For every $R>0$,
let $U_R$ be the maximal solution of
\begin{equation}\label{B-full}
 -\Gd v +g(v)=|f| \quad\text{in }B_R.
\end{equation}
By \rth{II}, $U_R$ is a large solution. Clearly $R\mapsto U_R$ decreases as $R$ increases. Therefore
$U=\lim_{R\tin}U_R$ is the maximal solution of
$$ -\Gd v +g(v)=|f| \quad\text{in }\BBR^N.$$
Similarly, let $W_R$ be the maximal solution of
\begin{equation}\label{B-full'}
-\Gd w +\tilde g(w)=|\tilde f|\quad\text{in }B_R,
\end{equation}
so that $W=\lim_{R\tin}W_R$ is the maximal solution of
$$ -\Gd v +\tilde g(v)=|\tilde f| \quad\text{in }\BBR^N.$$
If $u$ is any solution of (\ref{main-f}) in $B_R$ then $u\leq U_R$ and $\tilde u\leq W_R$ so that
\begin{equation}\label{B-ineq2}
  \tilde W_R\leq u\leq U_R.
\end{equation}
\note{Existence}
Let $k>0$, put $f_k=\min(|f|,k)\sign f$ and denote by $W_{k,R}$ and $U_{k,R}$ the maximal solutions defined above,
with $f$ replaced by $f_k$. Then $W_{k,R}$ and $U_{k,R}$ are locally bounded and increase with $k$.
Consequently, if $\{u^k_R:R>0\}$
is a family of functions such that $u^k_R$ is a solution of (\ref{main-f})  in $B_R$,
with $f$ replaced by $f_k$, this family is locally uniformly bounded.
This means  that, for every compact set $K$, there exits $R_k(K)>0$ such that
$\{u^k_R:R>R_k(K)\}$ is uniformly bounded in $K$.  Therefore there exists a sequence $R_j\tin$ such that
$\{u^k_{R_j}\}$ converges locally uniformly to a solution $u^k$ of (\ref{main-f})  in $\BBR^N$,
with $f$ replaced by $f_k$. By (\ref{B-ineq2}), the family of solutions $\{u^k:k>0\}$ is dominated
(in absolute value) by a function in $L^1_{loc}(\BBR^N)$ and it is non-decreasing.
Consequently $u=\lim u^k$ is a solution of (\ref{main-f}) in $\BBR^N$.\\[2mm]
\note{Uniqueness} Under the assumptions of (ii) $u\equiv 0$ is a solution of (\ref{main-0}) in $\BBR^N$
and it is easy to see that this is the only solution. Therefore the uniqueness statement for  (\ref{main-f})
follows  by the same argument as in the proof of \rth{IV}.
\eproof

\enddocument